# Digit Statistics of the First $\pi^e$ Trillion Decimal Digits of $\pi$

P. Trueb
DECTRIS Ltd., Baden-Dättwil, Switzerland,

*Abstract*—**The mathematical constant $\pi$ has recently been computed up to 22'459'157'718'361 decimal and 18'651'926'753'033 hexadecimal digits. As a simple check for the normality of $\pi$, the frequencies of all sequences with length one, two and three in the base 10 and base 16 representations are extracted. All evaluated frequencies are found to be consistent with the hypothesis of $\pi$ being a normal number in these bases.**

## I. Introduction

For a number to be normal in base $b$, every sequence of $k$ consecutive digits has to appear with a limiting frequency of $b^{-k}$ in the numbers' base-$b$ representation. It is supposed that $\pi$ is normal in any base, but a proof is still lacking [1]. As a consequence, new record computations of $\pi$ are often used to perform empirical consistency checks for the normality of $\pi$ [2]. Recently a new world record computation has been performed with the *y-cruncher* code [3]. This computation encompasses 22'459'157'718'361 decimal and 18'651'926'753'033 hexadecimal digits [4], adding about 70% more digits to the former record.

## II. Results

Figures 1 to 6 show the distributions of the frequencies of all sequences up to length 3 in the decimal and hexadecimal representations of $\pi$. The red and blue regions show the expected one and two standard deviation bands around the limiting frequency $b^{-k}$ assuming the occurrences of the sequences to follow a binomial distribution. All distributions show no significant irregularities. In particular there is no observed frequency deviating more than four standard deviations from the expected frequency. The expected and observed variances of the frequency distributions are listed in Table 1. All observed variances match the expected values, the maximum deviation amounting to 1.33 standard deviations.

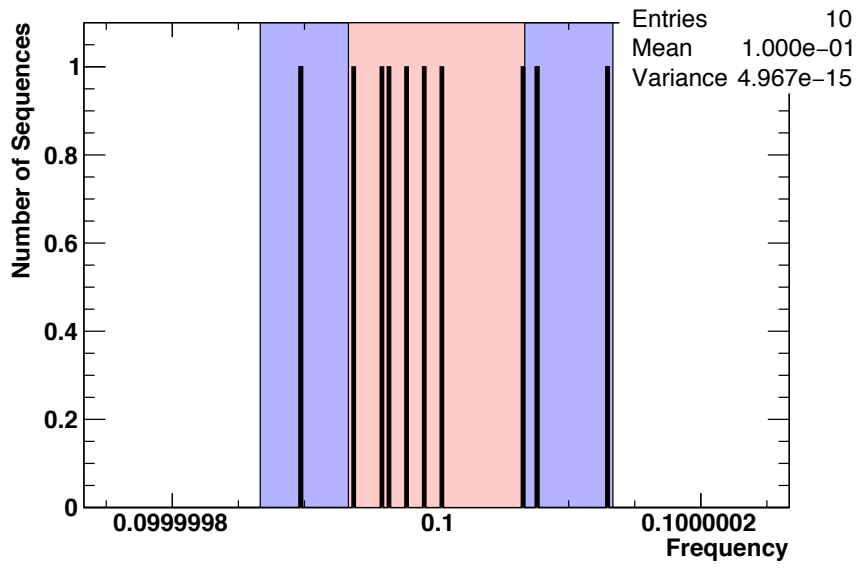

**Figure 1** Frequencies of the digits 0 to 9 in the decimal representation of π.

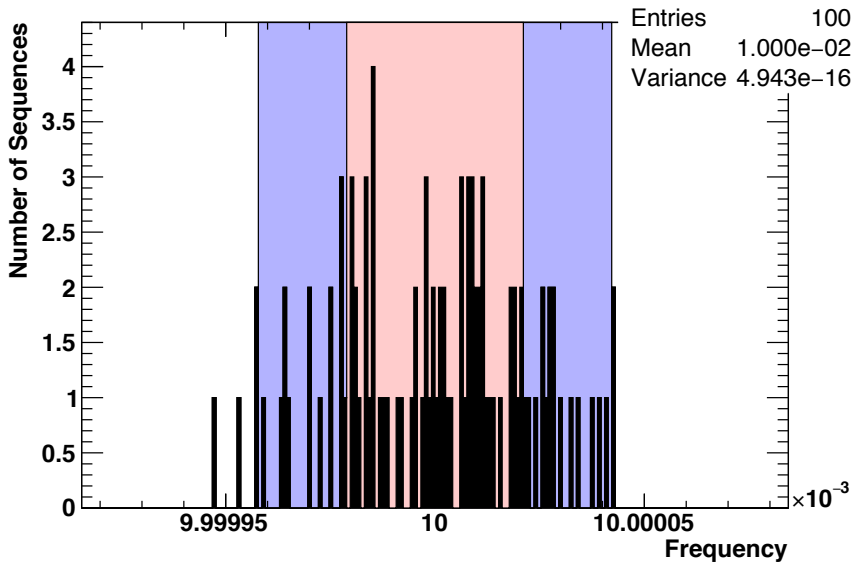

**Figure 2** Frequencies of all sequences of length 2 in the decimal representation of π.

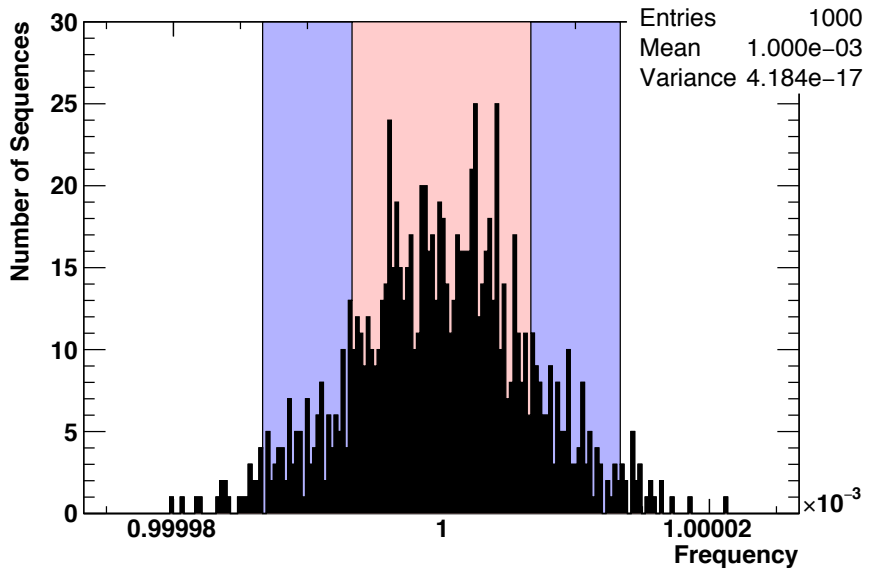

**Figure 3** Frequencies of all sequences of length 3 in the decimal representation of π.

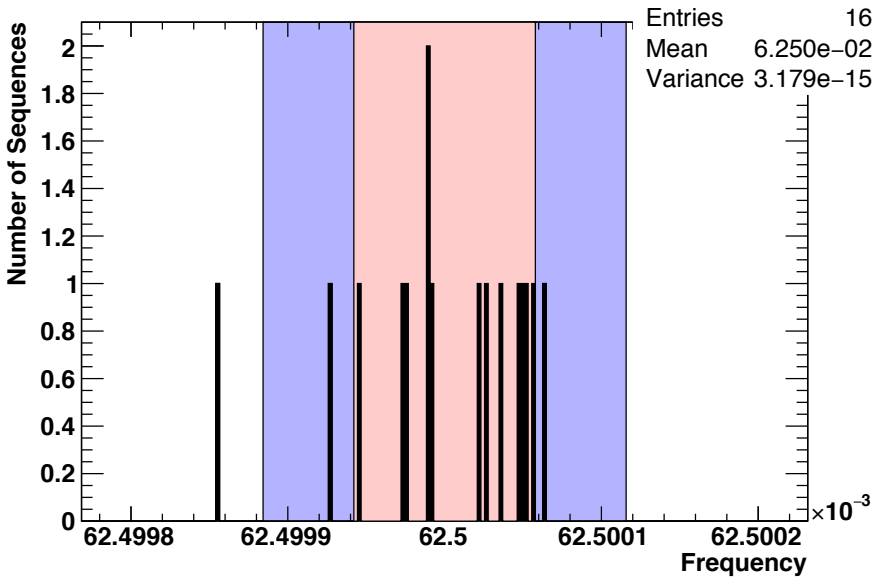

**Figure 4** Frequencies of the digits 0 to F in the hexadecimal representation of π.

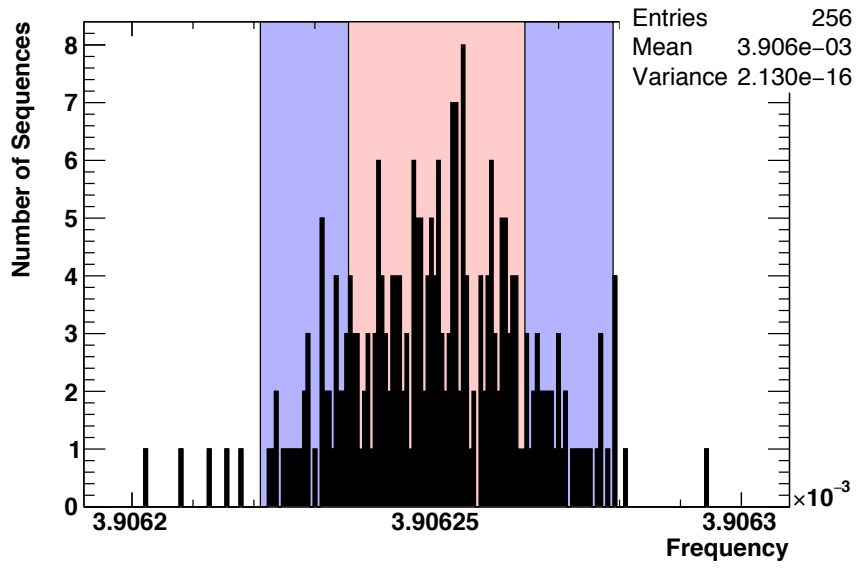

**Figure 5** Frequencies of all sequences of length 2 in the hexadecimal representation of π.

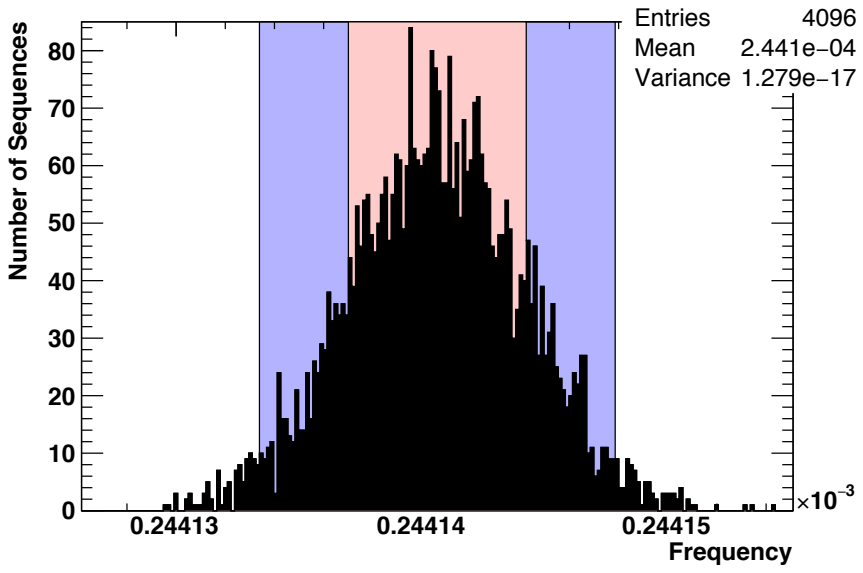

**Figure 6** Frequencies of all sequences of length 3 in the hexadecimal representation of π.

**Table 1** Expected and observed variances of the frequencies of all sequences up to length 3 in the decimal and hexadecimal representations of π.

| Base | Sequence Length | Expected Variance of Frequencies | Observed Variance of Frequencies | Deviation [σ] |
|---|---|---|---|---|
| 10 | 1 | $(4.0\pm1.9)\times10^{-15}$ | $5.0\times10^{-15}$ | -0.51 |
| 10 | 2 | $(4.41\pm0.63)\times10^{-16}$ | $4.94\times10^{-16}$ | -0.85 |
| 10 | 3 | $(4.45\pm0.20)\times10^{-17}$ | $4.18\times10^{-17}$ | 1.33 |
| 16 | 1 | $(3.1\pm1.1)\times10^{-15}$ | $3.2\times10^{-15}$ | -0.03 |
| 16 | 2 | $(2.09\pm0.18)\times10^{-16}$ | $2.13\times10^{-16}$ | -0.24 |
| 16 | 3 | $(1.309\pm0.029)\times10^{-17}$ | $1.279\times10^{-17}$ | 1.02 |

## III. Conclusions

The frequencies of all sequences up to length 3 in the first 22’459’157’718’361 decimal and 18’651’926’753’033 hexadecimal digits of π are found to be consistent with the hypothesis of π being a normal number in base 10 and base 16.